\documentclass[12pt]{article}
\usepackage{dsfont}
\usepackage{amsmath}
\usepackage{graphicx,graphics}
\usepackage{float}
\usepackage{subfigure}
\usepackage{bbm}
\usepackage{amsmath}
\usepackage{amsfonts}
\usepackage{latexsym}
\usepackage{amssymb}
\usepackage{verbatim}
\def\be{\begin{equation}}
\def\ee{\end{equation}}

\def\o{\over}
\def\ba{\begin{array}}
\def\ea{\end{array}}
\def\bea{\begin{eqnarray}}
\def\eea{\end{eqnarray}}
\def\h{\hbox}

\def\la{\langle}
\def\ra{\rangle}
\def\non{\nonumber}

\begin{document}
\begin{center}
{\large{\bf Balancing on the edge, the golden ratio, the  
		 Fibonacci sequence and their generalization }} \\
Gautam Dutta$^{1,4}$, Mitaxi Mehta$^2$, Praveen Pathak$^3$ \\
$^1${\it Dhirubhai Ambani Institute of Information and Communication Technology,\\
     Gandhinagar 382007,  India.}\\
$^2${\it SEAS, Ahmedabad University, Navrangpura, Ahmedabad 380009, India}    \\
$^3${\it Homi Bhabha Center for Science Education, Tata Institute 
	of Fundamental Research, Mankhurd, Mumbai, 400088, India}\\
$^4$ E-mail: gautam\_dutta@daiict.ac.in     
\end{center}

\begin{abstract} 
\noindent
 The golden ratio and Fibonacci numbers are found to  occur in 
 various aspects of nature.
 We discuss the occurrence of this ratio  in an interesting physical 
 problem concerning center of masses in two dimensions. The result is shown 
 to be independent of the particular shape of the object. The approach taken 
 extends naturally to higher dimensions. This leads to ratios similar to the golden 
 ratio and generalization of the Fibonacci sequence. 
 The  hierarchy of these ratios with 
 dimension and the
 limit as the dimension tends to infinity is discussed using the physical 
 problem. 
\end{abstract}

\noindent
\section{Introduction}
\begin{figure}[h]
  	\begin{center}
  	\includegraphics[scale=1.0]{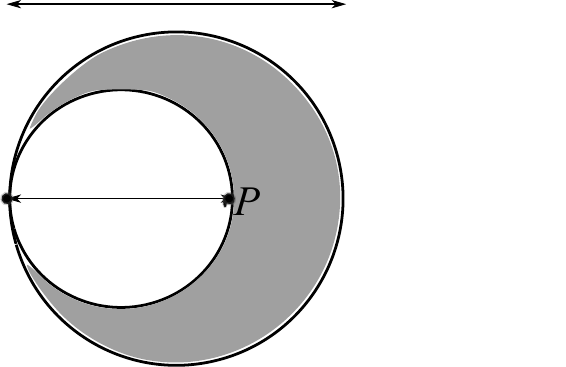}
  	\put(-130,50){$d'$}
  	\put(-115,106){$d$}
  	\put(-175,46){$O$}
  	\caption{A circular cavity  of diameter $d'$ is cut out from a 
  	circular disc of diameter $d$. The inner circle is tangential 
  	to the outer circle at the point $O$. The c.m of the  
  	crescent like remaining portion of the disc is at $P$.}
  	\label{circle}
  	\end{center}
  \end{figure}

In a basic physics course students handle a number of interesting problems on center of mass of laminar bodies.
In \cite{golden} the authors discuss an interesting connection of center of mass (c.m) of uniform laminar bodies to a very interesting irrational number, namely, the golden ratio. 
When a circular disk of diameter $d'$ is removed from a bigger circular disk of diameter $d$, touching the bigger circle
at a single point $O$,  it 
leaves a crescent like structure as shown in 
Fig. \ref{circle}. 
If we demand the c.m of the remaining portion of the disc to be  at $P$, a point exactly at 
the inner edge of the crescent, then 
the ratio of the diameter of the bigger to the smaller  circle is the golden ratio, i.e., 
${d\o d'}={1+\sqrt{5} \o 2}$. 
This result is shown and suitably extended to the case of an even sided regular polygon. 
Circles and regular polygons have rotational symmetry.  
In \cite{gautam} the result is extended to an ellipse which don't have rotational symmetry. 
In this article we generalize this for a uniform planar object of arbitrary shape. We also generalize the result beyond two dimensions whereby the said property of the c.m is related to irrational numbers similar to the golden ratio. 
The sequence of these irrational numbers converges to 2 as the dimension tends to infinity. \\

\noindent
A circular cavity is similar in shape to the circular disc 
from which it is cut since all circles are similar. 
Hence for a general non-circular planar object, we require
the cavity  that is cut to be  similar in shape to the 
original object, internally tangent and to have the same 
orientation as the original object.
This is described in Fig. \ref{general}. 
 Here, by similarity of 
 shapes we mean one is a scaled version of the other. 
We will consider the planar object to be convex. This ensures that the 
kind of cavity, as described,  lie completely  within 
the object. 
The analysis is applicable also to non convex objects where this 
condition is satisfied.\\

\noindent
The golden ratio \cite{dunlap} is the positive root of the quadratic polynomial 
$P_2(x)=x^2-x-1$. This irrational number, denoted as $\varphi$, 
mysteriously crops up in various natural phenomena \cite{mario}, 
geometrical  constructions \cite{steinhardt} and physical problems
\cite{srinivasan}. 
One of the most common geometrical description of the 
golden ratio is the ratio of length $l$ to the breadth $b$ of a 
rectangle such that we are left with a rectangle with the same ratio 
of length to breadth
after 
a square of side $b$ is removed from it, i.e, ${b\o l-b} = {l\o b}$. 
This ratio is also the asymptotic ratio of the $n^{th}$ and the $(n-1)^{th} $ term of a 
Fibonaccci sequence. A Fibonacci sequence is a sequence 
of numbers such that every term since the third term is a
sum of the previous two terms. 
The first two terms are arbitrary \cite{keith}. Usually one considers 
a Fibonacci sequence of positive integers, for e.g. 
$\{1, 1, 2, 3, 5, 8, 13, 21, 34,....\}$. The numbers in a Fibonacci sequence are called Fibonacci numbers. These numbers are found to occur in several places in nature such as the branching of stems in trees \cite{douady} and the spiral arrangement of buds in sunflower or the scales on pineapple \cite{brousseau}. \\

\noindent
We extend this analysis to three and higher dimensions and obtain 
higher degree analogues of the polynomial $ P_2(x)$. 
A connection is made with a kind of generalization of the Fibonacci 
sequence and that of the golden ratio.  \\

\noindent
\noindent
\section{Extension to an arbitrary shape in two dimension}
  \begin{figure}[h]
  	\begin{center}
  		\includegraphics[scale=1.0]{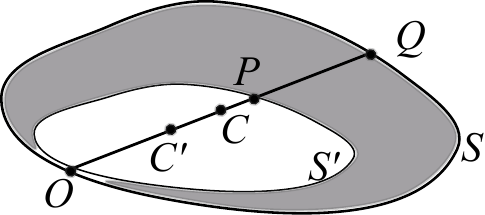}
  		\caption{A laminar object bounded by the curve $S$. $C$ is the c.m of the object. $S'$ is the edge of the cavity. $S'$ is similar to the curve $S$ and scaled around the common point $O$. $C'$ is the c.m of the cavity.   $P$ is the c.m of the remaining portion of the object lying on the inner edge $S'$. $Q$ is a point corresponding to $P$ on the outer edge $S$. The points $O,~C',~C,~P$ and $Q$ are collinear.   }
  		\label{general}
  	\end{center}
  \end{figure}
\noindent
Let $S$ be the boundary of the given object and $S'$ be the boundary of 
the cavity. They are closed curves and they touch each other at $O$. 
$S$ and $S'$ are similar, with $S'$ being the scaled down version of $S$. 
The common point $O$ is the center of scaling.   
$C$ is the c.m of the complete object. 
 Every point on the outer curve $S$ has a  corresponding point on the 
inner curve $S'$. 
Let $OQ$ be a chord of the closed curve $S$ that passes through the c.m 
$C$. 
$P$ is a point corresponding to $Q$ on $S'$. 
Let $C'$ be the c.m of the region removed by the cavity. 
Due to  similarity and identical orientation of $S$ and $S'$, 
$C'$ lies on the same chord $OQ$.\\
 
\noindent
Let $OQ= d$ and $OP= d'$. We will call $d$ and $d'$ as diameter.  
We make the following observations: The area of the given planar object
can be expressed as $\alpha d^2$ where $\alpha $ is a geometrical 
constant depending upon the shape of the object. For e.g in the case of 
a circle $\alpha= \pi/4$ when $OQ$ is the diameter. 
For a square $\alpha = 1$ when $OQ$ is the length of its side while 
$\alpha = 1/2$ when $OQ$ is the diagonal of a square. 
The area of the cut portion is 
$\alpha d'^2$ while the area of the shaded region is $\alpha (d^2-d'^2)$. Let $OC= \beta d$ with $0 < \beta < 1$.
Since the shape of the cavity is similar to that of  object, $OC' = \beta d'$. Like $\alpha,~~\beta$  is also a purely geometrical factor.
The center of mass of the shaded part will lie on the line joining 
$C$ and $C'$. Since we demand this center of mass to be on the edge of the 
cavity, it has to be the point $P$. 
Let $\sigma$ be the planar mass density.
Then the equation relating the three c.m is given as 
\be
\sigma \alpha(d^2-d'^2)d' + \sigma \alpha d'^2 \beta d'      \label{basic2d}
              = \sigma \alpha d^2 \beta d  
\ee
With $d=xd'$ this may be rewritten as   
\be
 (x^2-1)=\beta (x^3-1)
\ee
$x=1$ is an obvious solution of the above polynomial equation. This means the 
cavity is almost the size of the object. This can be understood only  as $x\to 1_+$. 
In this limit, the point $P\to Q$ and most of the mass lies between these points, 
hence near  $P$. This is not very interesting. Hence we look for other solutions. 
 Factoring
out $x-1$ from both sides gives the polynomial equation 
\be 
\beta x^2 + (\beta-1)x + (\beta -1) =0                   \label{gen}
\ee
This is the most general form in two  dimensional case. 
A few comments are in order. The roots of the quadratic polynomial 
in this equation are always real for $0 \le \beta \le 1$. 
The two roots are of opposite signs since the product of the roots 
is ${\beta-1\o \beta}$ which is negative. 
The absolute value of the positive root 
is greater than the absolute value of the negative root since the sum of the
roots is equal to ${1-\beta \o \beta}$, which is positive. These features will be 
observed even in higher dimensions.\\ 

\begin{figure}[h]
  	\begin{center}
  		\includegraphics[scale=1.0]{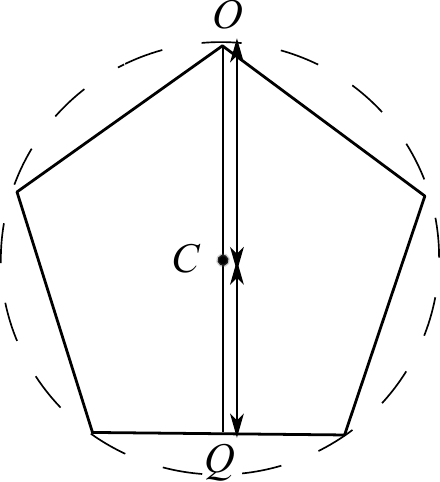}
  		\put(-55,80){$r$}
  		\put(-55,37){$r\cos{\pi \o n}$} 
  		\put(-69,-15){(a)} 		
  		\hspace{3cm}
  		\includegraphics{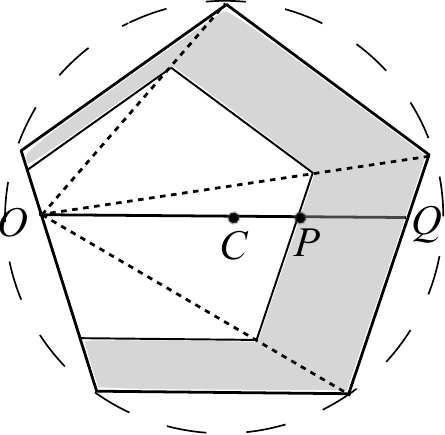}
  		\put(-69,-15){(b)}
  		\caption{An odd sided polygon. 
		(a): For the vertex excision distances are measured from the point $O$.  
		So $\beta= {1\o 1+\cos(\pi / n)}$. 
  		For the midpoint excision distances are measured from the point $Q$. So $\beta= {\cos(\pi / n)\o 1+\cos(\pi/ n)}$ .
  	    (b): $\beta=1/2$. Neither the vertices nor the mid points of the common sides of the two polygons coincide. The excission is assymetric. $P$ is the balancing point and ${OQ\o OP} = {1+\sqrt{5}\o 2}$.}
  		\label{polygon}
  	\end{center}
  \end{figure}

\noindent
Interestingly Eq. (\ref{gen}) is not dependent on the structure 
parameter $\alpha $. It only depends on $\beta $. For e.g, for circle,
even sided regular polygons and ellipse, the parameter $\alpha$ are 
different but the parameter $\beta =1/2$. Substituting this in the above equation gives 
\be 
x^2-x-1 =0                                        \label{symm}
\ee
for which the golden ratio is a positive root as mentioned above. These cases, discussed in \cite{golden} and \cite{gautam}, 
have a symmetry whereby the c.m $C$ is the midpoint of the chord $OQ$ independent of its orientation. 
All objects don't have such a symmetry and this is precisely one of the way we are generalizing here. 
If $C$ is not the midpoint of $OQ$ then $\beta_1 = {OC\o OQ} $ is either less than $1/2$ or greater than $1/2$, see Fig. \ref{general}. 
If the cavity is excised so that $Q$ is the common point then $\beta_2= {QC\o OQ} = 1-\beta_1$. 
As we rotate the chord $OQ$ about $C$, $\beta$ changes continuously from $\beta_1$ to $\beta_2$, 
thus taking the value $1/2$ at some orientation of $OQ$. 
So for any object we can find a chord along which Eq. (\ref{symm}) is satisfied. 
A case in point is an odd sided polygon (Fig. \ref{polygon}). 
Here, for the vertex excission, i.e, when the common point between the polygon and the excised cavity is the vertex $O$  in 
Fig. \ref{polygon}(a),  $\beta = {OC\o OQ}= {1\o 1+\cos(\pi /n)}$. 
For the case of mid point excission, i.e, when the common point
is the midpoint  $Q$ of a side of the polygon,  
$\beta = {QC\o OQ}={\cos(\pi /n)\o 1+\cos (\pi / n)}$. 
These cases are discussed in \cite{golden} whereby it is shown that the ratio ${d\o d'}$ required for balancing at the edge is not the golden ratio. 
As we rotate the chord $OQ$, $C$ becomes the midpoint of $OQ$ at a particular orientation. 
One such orientation is shown in Fig. \ref{polygon}(b).   
The way to make the excission is also shown.
$P$, a point on the inner edge, is the c.m of the  shaded region. 
The ratio, ${d\o d'}= {OQ\o OP}$ is the golden ratio in this case. 
Note that the excission is not symmetric, unlike the vertex and the midpoint excission, discussed in \cite{golden}. 

\noindent
Eq. (\ref{gen}) showed that the particular shape of the uniform laminar body doesn't matter. Now we know that for any shape we can find a chord through the c.m  along which 
$\beta = 1/2$ and  Eq. (\ref{gen}) takes the simple form, Eq. (\ref{symm}), whose positive root is ${1+\sqrt{5} \o 2}$. 
The golden ratio $\varphi$ is the asymptotic limit of the  ratio of successive terms of the 
Fibonacci sequence $F_n$ which is defined as $F_n=F_{n-1}+F_{n-2}$. 
So for large $n$, if $F_{n-2}=c$ then $F_{n-1}=c\varphi$ and $F_n=c\varphi^2$. Since these 
are the terms of the Fibonacci sequence we have 
\be c\varphi^2= c\varphi + c ~~~\h{i.e.},~~~\varphi^2-\varphi-1=0 \ee 
So the golden ratio satisfies Eq. (\ref{symm}). 
The other root,  $-1/\varphi$ is negative and hence not relevant to this problem. \\

\noindent
The problem, as presented  here, makes it easy to generalize beyond two dimensions. 
The independence of the shape and the existence of a chord along which $\beta = 1/2$ is easily extended to higher dimensions.
The visualization is same as given in Fig. \ref{general}. 

\section{Generalization to higher dimensions}

We start with modifying Eq. (\ref{basic2d}) for 
a homogeneous object in $k$ dimension. 
The $k$ dimensional volume of the object 
can be generally written as $\alpha d^k$ where $\alpha$ is again a purely geometric 
factor depending on the generic shape of the object and $d$ is the length scale in 
terms of which the volume can be determined. 
Likewise the volume of the cavity 
which will be removed is given by $\alpha d'^k $ where $d'$ is the corresponding 
length scale of the cavity. $OC=\beta d$ and $OC'=\beta d'$  where $0 < \beta < 1$
as discussed in the 2 dimensional case. So after removing the cavity if we want the 
c.m of the remaining object to be at $P$,  a point on the inner edge of the cavity then 
we get the following $k$ dimensional equivalent of Eq. (\ref{basic2d})
\be
\rho \alpha(d^k-d'^k)d' + \rho \alpha d'^k \beta d'      \label{basickd}
= \rho \alpha d^k \beta d  
\ee
where $\rho $ is the uniform $k$ dimensional volume density of the material.
Writing $d=xd'$, Eq. (\ref{basickd}) reduces to 
\be
 x^k-1 +\beta = \beta x^{k+1}                         
\ee
which can be rearranged as 
\be
 \beta (x^{k+1}-1)= x^k-1      
\ee
$x=1$ is a solution of the above equation which as discussed earlier is not of interest 
to us. 
So we divide both sides by $x-1$. 
We note that $1+x+....+x^{k-1} = {x^k-1\o x-1} $. 
Using this result on both sides we get 
\be
\beta x^k +(\beta-1)x^{k-1}+(\beta-1)x^{k-2}+....+(\beta-1)x+(\beta-1) =0
\ee
This is the $k$ dimensional equivalent of Eq. (\ref{gen}). Again we take the special 
case $\beta =1/2$. This gives us the equivalent of Eq. (\ref{symm}). 
\be
 x^k -x^{k-1} - x^{k-2}-....-x-1 =0              \label{symm-kd}
\ee
We had seen that Eq. (\ref{symm}) describes  the asymptotic behavior of the Fibonacci sequence. 
Likewise we discuss the relation
between Eq. (\ref{symm-kd}) and the generalized Fibonacci sequence.
  
\noindent
\section{Generalized Fibonacci sequence and its connection to higher 
	dimensional cases}  
There are various ways to generalize the Fibonacci sequence. We would be interested in a generalization to  a sequence of non-negative integers where the $n^{th}$ term 
is a sum of $k$ previous terms, i.e., $F_n=F_{n-1}+F_{n-2}+.....+F_{n-k}$. This is called the $k$- generalized Fibonacci sequence and the terms of this sequence are called the $k$-nacci numbers. 
So $k=2$ corresponds to Fibonacci numbers. $k=3$ gives the tribonacci numbers and likewise we have the tetranacci and pentanacci numbers.
Asymptotically for large $n$ this sequence also behaves like a geometric sequence 
with common ratio say $\zeta$~ \cite{bacani,tony}. So for large $n$ 
if $F_{n-k}=c$ then 
$F_{n-k+1}=c\zeta, F_{n-k+2}=c\zeta^2, ...., F_{n-1}=c\zeta^{k-1}$ and $F_n=c\zeta^k$. 
Substituting these in the generalized Fibonacci terms we get 
\be
\zeta^k=\zeta^{k-1}+\zeta^{k-2}+......+\zeta + 1
\ee
$\zeta$ satisfies Eq. (\ref{symm-kd}). It is the equivalent 
of the golden ratio in $k$ dimension and called the $k$-nacci constants.  

\noindent
Let us denote the polynomial on the l.h.s of Eq. (\ref{symm-kd}) as $P_k(x)$ which is a 
polynomial of degree $k$. 
Let us look at the kind of roots $P_k(x)$ can have.
$P_k(0) = -1 , ~~\forall~~ k $. For sufficiently large positive value 
of $x$, $P_k(x) > 0$. So each of these polynomials have a positive root. 
For even $k$, $P_k(x) >0$ also for sufficiently large negative value of 
$x$ . So these polynomials will also have a negative root which is not of relevance to us.
Can $P_k(x)$ have more than one positive root? In that case which of the positive roots will be the $k$-nacci constants? In fact we can show that it has exactly one positive root. 
Let $a$ be a positive root of $P_k(x)$. Then $a^k = 1+a+a^2+....+a^{k-1}$. 
It can be shown that  $x^k$ grows faster than $1+x+x^2+....+x^{k-1}$ for $x > a$. Hence there will be no other positive roots of $P_k(x)$.
  Hence for each $k$
we get exactly one positive  root, which are physically relevant to us. 
They are $k$ dimensional equivalent of the golden ratio, called the
$k$-nacci constants. We denote them as $\varphi_k$.
We will see that $\varphi_k$ is an irrational number between 1 and 2 for all $k$.
The coefficients of $P_k(x)$ are all integers and the leading coefficient (the coefficient of $x^k$) is 1. 
Such polynomials are called integer monic. They can be factored into integer monic \cite{int-poly}. 
So  the real roots of $P_k(x)$ are either integers or irrational. 
When $k=1,~P_1(x)=x-1$. So $\varphi_1 =1$. 
It can be shown algebraically that the positive roots of $P_k(x)$ lie between 1 and 2 but in the spirit of this article we give a strightforward justification using the generalized fibonacci sequence.  For $k > 1$ it is obvious that $F_n > F_{n-1}$  in a $k$-generalized Fibonacci sequence of positive numbers. 
Also  ${F_n < 2 F_{n-1}} $. This can be p{\large }roved as follows: 
\bea
F_{n-1} &=& F_{n-2}+F_{n-3}+....+F_{n-k}+F_{n-(k+1)} \non \\
F_{n} &=& F_{n-1}+(F_{n-2}+....+F_{n-k})  \non \\
      & < & F_{n-1}+F_{n-1} =2F_{n-1}    
\eea   
So  $1 < {F_n\o F_{n-1}} < 2$ for $k >1$. 
This implies $\varphi$ can not be an integer for $k > 1$. Hence they 
have to be an irrational number between 1 and 2. 
So we conclude that each of the $\varphi_k$ for $k >1$ is an irrational number between 1 and 2. 
$\varphi_1=1$. $\varphi_2= {1+\sqrt{5} \o 2}\approx 1.6180$ is the well known 
golden ratio, $\varphi_3 \approx 1.8393$ is the tribonacii constant, $\varphi_4 \approx 1.9276$ is the tetranacii constant and
 so on.  
$\varphi_k \to 2 $ as $k\to \infty$. 
This limit  is very interesting and simple to understand. 
Let us take an example of such a generalized Fibonacci sequence 
\be\la 0,~0,~0,~ 1,~ 1,~ 2,~ 4,~ 8,~ 16,~ .........\ra\ee  
In this sequence we can start with as many 0 as we like initially. 
Subsequent to the first non-zero entry, every term is the sum of all the 
previous terms. Asymptotically, this sequence behaves as the 
geometric sequence $2^n$. Hence the asymptotic ratio of successive terms is 2. \\

\noindent
The generalized golden ratios have a hierarchy, 
$\varphi_1 < \varphi_2 < \varphi_3 < ..... <2$.
This hierarchy can be understood algebraically but since these ratios 
occur in a physical problem, it is interesting to have a physical 
understanding for this.
 The mass of the object is proportional to $d^k$ 
and that of the cavity goes as  $d'^k$ in $k$ dimension.
The removal of the cavity shifts the c.m from $C$ to $P$
(See Fig. (\ref{general})). 
The fraction of the cavity mass
increasingly diminishes in comparison to the total mass 
with rising $k$ since $d' < d$. 
This would mean that the shift $C$ to $P$  decreases as $k$ increases. 
As $d'$ decreases with $k$ , $\varphi_k={d\o d'}$ increases with $k$. 
This explains the hierarchy of $\varphi_k$ with $k$. 
As $k\to \infty$, the mass of the cavity becomes negligible compared to the bulk.
Hence $P$ coincides  with $C$. In this case ${d\o d'}\approx {OQ\o OC} = {1\o \beta}$. 
If $\beta = {1\o 2} ,~~{d\o d'}\approx 2$. So $\varphi_k \to 2$ as $k\to \infty$.  

\section{Conclusion}
The physical problem of balancing a laminar object on its inner edge 
after suitably excising a self similar cavity within it, is found to be 
related to the golden ratio. We find that this phenomenon is independent 
of  the shape of the object thus emphasizing the role of this important 
irrational number in various aspects of nature.  The extension of this 
to a higher dimension $k$ gives a connection of this physical problem with the generalized Fibonacci sequence where every term is a sum of  
$k$ preceding terms. Due to this connection we give a physical 
explanation for the hierarchy of the generalized golden ratios 
$\varphi_k$ with $k$.   

\noindent
\begin{center}
{\bf ACKNOWLEDGEMENTS}
\end{center}
 This work is supported by the Physics olympiad programme, HBCSE-TIFR. 
 We thank Prof. Vijay A Singh for discussions and inspiring us with
 critical  questions during this work.

\end{document}